\documentclass[12pt,notitlepage]{amsart} 
\usepackage{amssymb,amsmath,epic,latexsym}

\newlength{\jmr}
\newlength{\baker}
\newlength{\jones}
\newlength{\sil}
\newlength{\foo}
\newtheorem{conj}{Conjecture}
\newtheorem{lemma}{Lemma}

\newtheorem{dfn}{Definition}

\newtheorem{mrt}{The MR Theorem}

\newtheorem{sie}{Siegel's Theorem}

\newtheorem{main}{Main Theorem}

\newtheorem{jst}{The JST Theorem} 

\newtheorem{thm}[main]{Theorem}
\newtheorem{rem}{Remark}	
\newtheorem{ex}{Example}

\newcommand{\cR}{\mathcal{G}}
\newcommand{\cT}{\mathcal{T}}
\newcommand{\tC}{\tilde{C}}

\newcommand{\supp}{\mathrm{Supp}}

\newcommand{\thth}{^{\underline{\mathrm{th}}}}

\newcommand{\Pro}{\mathbb{P}}
\newcommand{\Q}{{\mathbb{Q}}}
\newcommand{\bQ}{\overline{\Q}}
\newcommand{\R}{\mathbb{R}}
\newcommand{\C}{\mathbb{C}}
\newcommand{\biggy}{\mathbf{Big}}
\newcommand{\card}{\mathbf{ExactCard}}
\newcommand{\N}{\mathbb{N}}
\newcommand{\Z}{\mathbb{Z}}
 
\newcommand{\fii}{\varphi}

\newcommand{\codim}{\mathrm{codim}}

\newcommand{\newt}{\mathrm{Newt}} 
\newcommand{\htp}{\mathcal{HTP}}  
\newcommand{\seq}{\mathrm{MultiSub}}

\newcommand{\Zn}{\Z^n}

\newcommand{\Rn}{\R^n}
\newcommand{\Cn}{\C^n}

\renewcommand{\qed}{$\blacksquare$}

\newcommand{\cG}{\mathcal{G}}

\begin{document}

\title[Uncomputably Large Integral Points?]{Uncomputably Large Integral Points 
on Algebraic Plane Curves?$^\star$ }

\thanks{$^\star$ In: {\it Theoretical Computer Science}, special issue in 
honor of Professor Manuel Blum's 60$\thth$ birthday, to appear. } 


\author{J. Maurice Rojas}

\thanks{Partially supported by a National Science 
Foundation Mathematical Sciences Postdoctoral Fellowship and a 
Hong Kong CERG Grant.} 

\address{Department of Mathematics, City University of Hong Kong, 
83 Tat Chee Avenue, Kowloon, HONG KONG}

\email{mamrojas@math.cityu.edu.hk}

\dedicatory{Dedicated to Professor Manuel Blum on the 
occassion of his $60\thth$ birthday.} 

\date{\today}

\begin{abstract} 
We show that the decidability of an amplification of Hilbert's Tenth Problem 
in three variables implies the existence of uncomputably large 
integral points on certain algebraic curves.
We obtain this as a corollary of a new positive complexity result: 
the Diophantine prefixes $\exists\forall\exists$ and 
$\exists\exists\forall\exists$ are {\bf generically} decidable. This means, 
taking the former prefix as an example, that we give a 
precise geometric classification of those 
polynomials $f\!\in\!\Z[v,x,y]$ for which the question 
\[ \exists v\!\in\!\N \mathrm{ \ 
such \ that \ } 
\forall x\!\in\!\N \ \ \exists y\!\in\!\N \mathrm{ \ with \ } 
f(v,x,y)\!=\!0? \hspace{2cm} \mbox{} \] 
may be undecidable, and we show that this set of polynomials is quite small in 
a rigourous sense. (The decidability of $\exists\forall\exists$ was 
previously an open question.) The analogous result for the prefix 
$\exists\exists\forall\exists$ is even stronger. We thus obtain a 
connection between the decidability of certain Diophantine problems, height 
bounds for points on curves, and the geometry of certain complex surfaces and 
$3$-folds.  
\end{abstract} 

\mbox{}\\
\vspace{-1in}
\maketitle

\section{Introduction}
\label{sec:intro}

We derive new complexity-theoretic limits 
on what can be discerned about the set of 
integral points of a variety of low dimension. 
In particular, we exhibit a new family 
of decidable Diophantine sentences related 
to the remaining open cases of Hilbert's Tenth 
Problem. As a corollary, we obtain a Diophantine 
problem whose decidability implies the 
following surprising assertion: for a general algebraic plane curve 
$\{ (x,y)\!\in\!\C \; | \; f(x,y)\!=\!0\}$, it is impossible 
to express the size of the largest positive integral point 
as a Turing computable function of the degree 
and coefficient sizes of $f$. 

Finding such bounds is literally one of the holy grails of number theory.
Huge, but nevertheless computable, upper bounds have already been found for 
a number of important classes of curves, such 
as curves of genus one \cite{bakercoates}, Thue curves 
\cite{bakert}, hyperelliptic curves \cite{bakerh}, 
superelliptic curves \cite{brindza}, and 
certain rational curves \cite{poulaki}.\footnote{By a {\bf curve} (resp.\ 
{\bf surface} or {\bf $\boldsymbol{3}$-fold}) we will always mean a 
one-dimensional (resp.\ two-dimensional or three-dimensional) 
{\bf complex} zero set of some set of polynomial equations.} 
For example, it is known that for any polynomial equation of the form 
\[ y^2=a+bx+cx^2+dx^3, \] 
where $a,b,c,d\!\in\!\Z$ and $a+bx+cx^2+dx^3$ has three 
distinct complex roots, all {\bf integral} solutions must  
satisfy 
\[|x|,|y|\leq \exp((10^6H)^{10^6}),\] 
where $H$ is any upper bound on $|a|,|b|,|c|,|d|$ \cite{bakertran}.\footnote{ 
For any integral point $(x,y)\in\Z^2$, the quantity 
$\max\{|x|,|y|\}$ is usually called the {\bf height} of $(x,y)$. Also, the 
preceding height bound has since been considerably improved, 
e.g., \cite{schmidt}.} 

However, finding such bounds, even ones 
monstrously larger than those already known, for {\bf general} algebraic curves has been out of reach for decades.
Furthermore, the analogous question for algebraic surfaces, even 
in $\C^3$, has so far been addressed only through deep conjectures 
of Lang and Vojta \cite{lang,vojta}. 

Our first main theorem relates the decidability of certain 
Diophantine sentences in four variables with the computability of upper bounds on the size of integral points on algebraic curves. 
So let us briefly recall Hilbert's Tenth Problem in $n$ variables:
\vspace{-.5cm}
\begin{center}
{\bf 
\mbox{\hspace{-.2cm}``Decide whether an arbitrary 
$\boldsymbol{f\!\in\!\Z[x_1,\ldots,x_n]}$ has an 
integral root or not.''}}
\end{center} 
We will denote this well-known Diophantine problem by 
$\boldsymbol{\htp_\Z(n)}$. 
Similarly, the analogous problem where we wish to determine the existence of 
an integral root, with all coordinates {\bf positive}, will be denoted by 
$\boldsymbol{\htp_\N(n)}$. We will also need the following 
closely related functions.  
\begin{dfn}
\label{dfn:big}
For any subset $R\!\subseteq\!\C$ closed under addition 
and multiplication, define the functions  
\[ \biggy_{R,n},\card_{R,n} : \Z[x_1,\ldots,x_n] \longrightarrow \N\cup\{0,\infty\} \] 
as follows: Let $S_f$ be the hypersurface $\{ (x_1,\ldots,x_n)\!\in\!\C^n \; | 
\; f(x_1,\ldots,x_n)\!=\!0 \}$. Then $\biggy_{R,n}(f)$ is the 
supremum of $\max\{|r_1|,\ldots,|r_n|\}$ as $(r_1,\ldots,r_n)$ 
ranges over $\{(0,\ldots,0)\}$ and the set of all points of $S_f$ 
in $R^n$. Finally, $\card_{R,n}(f)$ is the number of points of $S_f$ 
in $R^n$.  
\end{dfn}
\noindent 
It is not hard to see that the decidability of $\htp_\N(n)$ implies the 
decidability of $\htp_\Z(n)$, so $\htp_\N(n)$ is at least 
as hard as $\htp_\Z(n)$. Similarly, the computability of 
$\biggy_{\N,n}$ (resp.\ $\card_{\N,n}$) easily implies 
the computability of $\biggy_{\Z,n}$ (resp.\ $\card_{\Z,n}$). 
Also, via brute-force enumeration, it is easy to see 
that $\card_{\N,n}$ (resp.\ $\card_{\Z,n}$) is computable 
iff $\biggy_{\N,n}$ (resp.\ $\biggy_{\Z,n}$) is computable.  
However, we also have the following more subtle fact.  
\begin{main}
At least one of the following two statements is {\bf false}:
\begin{enumerate}
\item{The function $\biggy_{\N,2}$ is Turing computable.}
\item{The Diophantine sentence  
\[ \exists u\!\in\!\N \ \ \exists v\!\in\!\N \ \ \forall 
x\!\in\!\N \ \ 
\exists y\!\in\!\N \mathrm{ \ with \ } f(u,v,x,y)\!=\!0 \] 
is decidable in the special case where the underlying 
$3$-fold $S_f$ contains a surface which is a bundle of 
curves (each with a genus zero component) fibered over a curve $C$ 
in the $(u,v)$-plane, where $C$ has infinitely many positive integral 
points. } 
\end{enumerate}
In particular, $\htp_\N(3)$ is a special case of the problem 
mentioned in statement (2). 
\end{main}
\noindent 
The geometric notions mentioned above are clarified in section 2. 
Alternative classes of $C$ for which Main Theorem 1 remains true 
are mentioned in sections 1.1 and 5.    

We note that Alan Baker has conjectured \cite[Section 5]{jones81} that 
$\htp_\Z(2)$ is decidable. Thus the truth of statement 
(1) above would imply an algorithm for deciding $\htp_\N(2)$, and thus a 
positive answer to Baker's conjecture as well. We also point out that the 
computability of $\card_{\Z,2}$ (and $\card_{\N,2}$) is {\bf still} 
an open question, in spite of the fact that explicit 
(albeit huge) upper bounds on $\card_{\Z,2}$ are known 
in many cases \cite{bomb,poulaki}. 
 
On the other hand, Z.\ W.\ Sun has proved that $\htp_\Z(11)$ is undecidable 
\cite{sunphd}. Also, Y.\  V.\ Matiyasevich 
has shown (the proof appearing in a paper of J.\ P.\ Jones \cite{jones9}) that 
$\htp_\N(9)$ is undecidable. However, although the decidability of 
$\htp_\N(1)$ and $\htp_\Z(1)$ is a simple algebraic exercise, 
the remaining cases of $\htp_\N(n)$ and $\htp_\Z(n)$, as of 
mid-1998, are still completely open.

Our result above thus tells us something new about the next harder 
(and open) cases of $\htp_\N(n)$. 

\begin{rem}
The computability of $\biggy_{\Z,n}$ (resp.\ 
$\card_{\Z,n}$) does {\bf not} trivially imply the computability 
of $\biggy_{\N,n}$ (resp.\ $\card_{\N,n}$): It is possible for 
$\biggy_{\Z,n}(f)$ (resp.\ $\card_{\Z,n}$) to be infinite and thus give 
us no decisive information about the value of $\biggy_{\N,n}(f)$ 
(resp.\ $\card_{\N,n}$). 
\end{rem}

This connection between height bounds and Hilbert's Tenth Problem 
points to an unusual possibility: The search for general effective 
height bounds for integral points on algebraic curves may be 
futile. Indeed, it would have perhaps been more interesting  
to prove the statement ``$\htp_\N(3)$ is decidable 
$\Longrightarrow \biggy_{\N,2}$ is uncomputable,'' or better 
still, ``$\htp_\Z(3)$ is decidable $\Longrightarrow \biggy_{\Z,2}$ 
is uncomputable.'' 
However, Main Theorem 1 is at least a first step in this direction.  We will 
comment further on strengthening Main Theorem 1 in the conclusion of this 
paper. 

While our first main result is negative in the sense that it implies 
undecidability for certain Diophantine sentences, its proof 
follows easily from our derivation of two positive results on Diophantine 
sentences. To describe these results, let us introduce the following notation: 
We say that ``the {\bf Diophantine prefix} $\exists v \forall x\exists y$ 
is decidable'' iff there is a Turing machine algorithm which decides 
the sentence 
\[ \exists v \forall x\exists y\mathrm{ \ with \ } f(v,x,y)\!=\!0 \]
for arbitrary input $f\!\in\!\Z[v,x,y]$, and where the quantification 
is over the {\bf \underline{positive}} integers. This notation 
extends in an obvious way to other combinations of quantifiers 
and variables such as $\exists v\exists y$, $\exists u\exists v\exists y$, 
etc. Finally, by {\bf generic} decidability, 
we will mean that a prefix is decidable when the input is 
restricted to an a priori fixed ``large'' set. This 
is made more precise below and in section 2. 

We will prove the following result. 
\begin{main}
The prefix $\exists v\forall x\exists y$ is {\bf generically} decidable. 
More precisely, it is decidable on the collection of those $f$ for which the 
underlying complex surface $S_f$ does {\bf not} have an irreducible 
component which is a bundle of curves (each with a genus zero 
component) fibered over the $v$-axis. 
\end{main} 
\noindent 
By simply considering those polynomials in $\Z[v,y]$, 
note that the prefix $\exists v\exists y$ (or, equivalently, 
$\htp_\N(2)$) is a special case of the prefix $\exists v \forall x\exists y$. 
It is also easy to see (cf.\ section 2) that the set of ``hard'' $f$ 
{\bf omitted} by our result above happens to include $\Z[v,y]$. 
Furthermore, via theorem \ref{thm:locusg} at the end of this section, 
we can algorithmically determine whether $f$ satisfies the above hypothesis. 

It should of course be pointed out that 
the decidability of $\exists\forall\exists$ was a completely 
open problem. In fact, J. P. Jones \cite{jones81} has conjectured 
that the prefixes $\exists\forall\exists$ and $\exists\exists$ 
are equivalent. Put another way, this is the conjecture that 
$\exists\forall\exists$ is decidable $\Longleftrightarrow 
\htp_\N(2)$ is decidable. So while we still haven't resolved the decidability of $\exists\forall\exists$, 
we now at least know a geometric characterization of where any 
potential obstruction to decidability may lie.  
In particular, it follows from a fundamental 
result of algebraic geometry that our hypothesis rules 
out certain {\bf ruled} surfaces, i.e., surfaces 
which are traced out by an infinite family of lines. 
The latter statement is also clarified in section 2. 

Our final main theorem is a seemingly paradoxical 
extension of the preceding result.
\begin{main}
The prefix $\exists u\exists v\forall x\exists y$ 
is generically decidable. More 
precisely, it is decidable on the collection of $f$ for which the underlying 
$3$-fold $S_f$ contains {\bf no} surface which is a 
bundle of curves (each with a genus zero component) fibered over 
a curve in the $(u,v)$-plane. 
\end{main} 
\noindent
We can algorithmically determine whether 
$f$ satisfies the preceding hypothesis as well, via theorem \ref{thm:locusg} at the end of this section. 

The ``near paradox'' arises from the following result of 
Y. V. Matiyasevich and Julia Robinson.\footnote{The paper \cite{jones81} 
contains many important results related to the MR Theorem, and for non-Russian 
readers may be a better reference than the original reference 
\cite{matrob}.}
\begin{mrt} 
\cite{matrob}
The quantifier prefix $\exists\exists\forall\exists$ 
is {\bf un}decidable, i.e., there is no Turing machine 
which decides for an {\bf arbitrary} input $f\!\in\!\Z[u,v,x,y]$ 
whether there is a $(u,v)\!\in\!\N^2$ such that 
$\forall x \; \exists y$ with $f(u,v,x,y)\!=\!0$. \qed 
\end{mrt} 
Thus, our generic decidability result is stronger for 
the prefix $\exists\exists\forall\exists$: We obtain a  
necessary geometric condition classifying those $f$ for which 
the above Diophantine sentence is undecidable. 
It is easy to see (cf.\ section 2) that this set of ``hard'' $f$ includes the 
prefix $\exists u \exists v \exists y$ (i.e., the problem 
$\htp_\N(3)$). However, this does {\bf not} 
necessarily imply that the prefix $\exists u \exists v \exists y$ 
is undecidable --- the undecidability of $\exists u \exists v \forall 
x \exists y$ may be due to {\bf other} polynomials in our exceptional locus. 
We also emphasize that the set of exceptional $f$ in Main Theorem 3 
is strictly larger than the set of $f$ considered in 
Main Theorem 1. 

The proofs of Main Theorems 2 and 3 are not difficult conceptually, 
but rely upon results of C.\ Runge \cite{runge,ayad}, C.\ L.\ Siegel 
\cite{siegel}, and A.\ Schinzel \cite{schinzel} on the distribution of integral 
points on curves. The necessary results are stated in 
section 1.1. The application of these results then relies on combining a 
geometric construction with a new, more effective 
characterization of genus zero for algebraic curves. The 
following definition makes this more precise.  

\begin{dfn}
\label{dfn:g0} 
Suppose $g\!\in\!\C[a_1,\ldots,a_m,x_1,\ldots,x_n]$ and 
let $a\!:=\!(a_1,\ldots,a_m)$. If we 
choose constants in $\C$ for all the $a_i$, we denote the 
corresponding specialization of $g$ to a 
polynomial in $\C[x_1,\ldots,x_n]$ by $g_a$. For any 
$g\!\in\!\C[a_1,\ldots,a_m,x_1,x_2]$, we then define the {\bf 
genus zero locus of $g$}, $\cR_g$, 
to be the set of all $a\!\in\!\C^m$ 
such that $S_{g_a}$ has an irreducible component with 
geometric genus zero.\footnote{ 
For convenience, we will sometimes use $x$ and $y$ in place 
of $x_1$ and $x_2$ in the bivariate case.}  
\end{dfn} 

By a celebrated theorem of Siegel (cf.\ section 1.1), 
detecting genus zero for a given curve (in many cases) is 
equivalent to detecting the existence of infinitely many integral
points. So the following theorem, proved in 
section 3, may be of independent interest. 
\begin{thm}
\label{thm:locusg} 
For any $g\!\in\!\bQ[a_1,\ldots,a_m,x,y]$, the locus $\cR_g$ is 
a quasi-affine variety, and the equations (and inequations) defining 
$\cR_g$ can be constructed effectively, e.g., by a Turing machine.  
\end{thm}
\noindent 
For example, a special case which is easy to derive from 
the basic theory of elliptic curves \cite{sil} 
is the following: If $g\!:=\!a_1y^2+a_2+a_3x+a_4x^3$, 
then the zero set of $g$ in $\C^2$ has an irreducible component 
of genus zero iff 
\[ a_1a_4(4a^3_3+27a^2_2a_4)\!=\!0.  \] 
Oddly, while there are certainly algorithms for computing 
the genus of the zero set of a given irreducible 
$f\!\in\!\bQ[x,y]$ (e.g., \cite{hoeij}), the effective geometric 
characterization of genus zero above appears to be new. 
So we present a proof of theorem \ref{thm:locusg} in section 3.  

In closing this first half of our introduction, we point out that 
our main theorems suggest that there is a deep connection between 
complex geometry and Diophantine complexity 
which has yet to be explored. In particular, 
we clearly need more refined geometric invariants 
to explicitly classify those curves (and surfaces) where we 
can hope to effectively study integral points. 

Main Theorems 2 and 3 are proved in section 4, and  
Main Theorem 1 is then proved in section 5. Some interesting 
open questions are briefly discussed in section 6. We now describe 
our necessary results on integral points more precisely.

\subsection{Curves with Many Integral Points} 
\label{sec:arse} \mbox{}\\
In this subsection, we will let $f$ denote a polynomial 
in $\Z[x,y]$. Let us quote Siegel's classification 
of curves having infinitely many integral points. 
\begin{sie}
\cite{siegel,lang,ayad,wow}
Let $C$ be a curve defined over $\Z$ and irreducible over $\C$. Then 
$C$ has infinitely many integral points $\Longrightarrow C$ has genus 
zero and at most two distinct points in $\Pro^2_{\bar{\Q}}\!\setminus\!\C^2$. 
Furthermore, we have the following partial converse: Any $C$ with 
genus zero and at most two distinct points at infinity will 
have infinitely many integral points in a sufficiently large finite algebraic 
extension of $\Z$. 
\qed 
\end{sie}
\begin{rem}
By the genus of a variety $V$ we will always mean the geometric genus of a 
smooth projective model for $V$. 
(Since geometric genus is a birational invariant, 
it will be independent of the chosen model.) The genus of a curve is described 
very nicely in \cite{miranda}, and the genera of higher dimensional varieties 
is defined in \cite{hart,kho78}. 
\end{rem}
\begin{rem} 
\label{rem:infinity}
For $C\!=\!S_f$, our definition of $C$ having a ``point at infinity''  
is simply that the compactified zero set of $C$ intersect 
$\Pro^2_{\bar{\Q}}\!\setminus\!\C^2$.  
So the ``points at infinity'' condition in Siegel's Theorem can 
actually be checked algorithmically, simply by considering the zero set of the 
homogeneous polynomial $t^{\deg f}f(x/t,y/t)|_{t=0}$ in 
$\bQ^2\!\setminus\!(0,0)$. 
\end{rem} 

Note that an immediate corollary of Siegel's Theorem, and our 
preceding remarks, is that the condition on $C$ in Main Theorem 1 can 
be replaced by the following:
\begin{center}
``...$C$ has a component of genus zero with at most two 
distinct points at infinity.''
\end{center} 
This version can be checked algorithmically, but gives a 
slightly larger exceptional case of $\exists\exists\forall\exists$ 
than the original condition. 
 
The final result on integral points we quote will allow us to 
efficiently decide a small, but highly non-trivial quantifier prefix. 
\begin{jst}
\cite{jones81,schinzel,tungcomplex}
The quantifier prefix $\forall\exists$ is 
decidable in polynomial time. More explicitly, 
given $P\!\in\!\Z[x,y]$, we have that 
$\forall x\; \exists y\; P(x,y)\!=\!0$ 
iff all of the following conditions hold: 
\begin{enumerate}
\item{The polynomial $P$ factors into the form 
$P_0(x,y)\prod^k_{i=1}(y-P_i(x))$ where $P_0(x,y)\!\in\!\Q[x,y]$ 
has {\bf no} zeroes in the ring $\Q[x]$, and for all $i$, 
$P_i\!\in\!\Q[x]$ and the leading coefficient of $P_i$ is positive.}  
\item{$\forall x\!\in\!\{1,\ldots,x_0\} \; \exists 
y\!\in\!\N$ such that $P(x,y)\!=\!0$, where $x_0\!=\!\max\{s_1,\ldots,s_k\}$,  
and for all $i$, $s_i$ is the sum of the squares of the coefficients 
of $P_i$.} 
\item{Let $d$ be the least positive integer such that 
$dP_1,\ldots,dP_k\!\in\!\Z[x]$ and set $Q_i\!:=\!dP_i$ for all $i$. Then the 
{\bf union} of the solutions of the following $k$ congruences 
\begin{eqnarray*}
Q_1(x) & \!\equiv & 0 \mod d \\
 & \vdots & \\
Q_k(x) & \!\equiv & 0 \mod d
\end{eqnarray*} 
is {\bf all} of $\Z/d\Z$. }
\end{enumerate} 
In particular, the above conditions can be checked 
in time polynomial in $\log(d)$, the heights of the 
coefficients, and the degree of $P$, 
via fast factorization of polynomials in $\leq\!2$ variables 
over $\Q$ and $\Z/d\Z$ \cite{cohen}. \qed 
\end{jst} 
\begin{rem}
The JST Theorem can be strengthened slightly in the 
following way: one can replace $d$ in condition (3) 
with {\bf any} positive integer $d'$ such that $d'P_1,\ldots,
d'P_k\!\in\!\Z[x]$. 
\end{rem} 

\section{Geometric Background}  
\label{sec:hurwitz}

We first point out that a complete account 
of computability, decidability, and Turing machines can be found in 
\cite{gj,hilbert10,bcss}. Also, our notion of ``input'' 
will be fairly standard: either the sparse encoding of 
polynomials (over $\Z$ or $\Q$) or the bit-wise encoding 
of algebraic numbers, for BSS machines over $\Z/2\Z$ \cite{bcss}. 
Finally, for most of the basic facts we will use from 
algebraic geometry, we refer the reader to 
\cite{hart,mumford,beau}. However, for the convenience of the reader, 
we will restate a few of the most central notions. 

\begin{rem} 
Throughout most of this paper, ``effectively computable'' and 
``algorithmic'' will be taken to mean Turing computable. 
\end{rem} 

Returning to the concept of sparse encoding, the following 
notation will be useful: For any $e\!\in\!\Zn$, we 
let $x^e$ denote the monomial term $x^{e_1}_1x^{e_2}_2\cdots 
x^{e_n}_n$. For any polynomial $f(x)\!=\!\sum_{e\in\Zn} 
c_e(a) x^e\!\in\!\C[a_1,\ldots,a_m,x_1,\ldots,x_n]$, we then let the {\bf 
support} of $f$, $\supp(f)\!\subset\!\Zn$, be the set of exponents
$\{e\!\in\!\Zn \; | \; c_e(a)\!\not\equiv\!0\}$. (We are implicitly 
considering the $x_i$ as variables to be solved for, and the $a_j$ as 
parameters we are free to choose.) Also, we will let the {\bf Newton 
polytope} of $f$, $\newt(f)\!\subset\!\Rn$, be the convex hull of (i.e., 
the smallest convex set containing) $\supp(f)$. 

An {\bf affine variety} is simply the 
complex zero set of a system of polynomial equations.
(So our varieties will not necessarily be 
{\bf reduced} or {\bf irreducible} \cite{hart,mumford}.) 
More generally, a {\bf quasi-affine} variety is the 
set of complex points satisfying any finite Boolean combination 
of polynomial equations and inequations. (Note that 
we mean $\neq$, not $<$ or $\leq$, when we say inequation.) 
In particular, when we say a set of indeterminates 
$\{a_1,\ldots,a_m\}$ is chosen {\bf generically}, 
we will mean that $(a_1,\ldots,a_m)\!\in\!\C^m\!\setminus\!W$ for  
some {\bf a priori fixed} quasi-affine algebraic 
subvariety $W$ (depending only on the property in question) 
with $\codim W\!\geq\!1$.
For most purposes, assuming a property holds  
generically also implies that the property occurs  
with probability $1$. (A classification of a 
broad class of probability measures on $\C^m$ 
for which this is true is not hard to derive.) 
We will also use the term {\bf variety} collectively 
for affine and quasi-affine varieties. 

As for the geometric language of our main theorems, 
let us recall the following definitions: A {\bf 
morphism} is simply a well-defined map from one variety to 
another, given by rational functions. When 
we relax the ``well-defined'' stipulation and allow 
our map to be undefined on a subvariety of codimension 
$\geq\!1$, we 
then obtain a {\bf rational} map. 
Also, a {\bf birational} map is a rational map with an 
inverse which is again a rational map. 

The inverse image of a point, for any 
given morphism, is usually called a {\bf fiber}. 
For any curve $C$, we then say that 
a variety $V$ is a {\bf bundle of curves 
fibered over $C$} iff there is a morphism 
$\varphi : V \longrightarrow C$ such that 
every fiber of $\varphi$ is a (not necessarily irreducible) curve.  
\begin{dfn} \cite{beau} 
Assume temporarily that all varieties are irreducible, 
nonsingular, and compact. Let $C$ be a curve. 
We then call a surface $S$ {\bf 
ruled over $C$} iff there is a morphism $\varphi : S \longrightarrow 
C$ with every fiber isomorphic to $\Pro^1_\C$. Similarly, 
we will call a surface $S$ {\bf ruled} iff there is a curve $C$ for which 
$S$ is ruled over $C$. Finally, an arbitrary algebraic surface $S$ 
(not necessarily irreducible, nonsingular, or compact) is said to be {\bf rational} iff $S$ is birationally 
equivalent to $\Pro^2_\C$. 
\end{dfn} 

The following two related facts will prove useful. 
\begin{thm} 
\label{thm:kho} 
\cite{kho78}\footnote{The version of theorem \ref{thm:kho} stated in 
\cite{kho78} is actually a bit different, but easily implies 
our version here via an application of theorem \ref{thm:all} below.} 
Suppose $f\!\in\!\C[a_1,\ldots,a_m,x_1,\ldots,x_n]$. Then, for 
generic $a$, the genus of $S_{f_a}\!\subseteq\!\Cn$ is exactly  
the number of lattice points in the interior of $\newt(f)$. \qed 
\end{thm}
\begin{thm}
\label{thm:rule}
If $S$ is a ruled surface, then its genus is zero.
Also, if $S$ is a bundle of curves (each with a genus zero 
component) fibered over another 
curve, then $S$ has a component which is birationally equivalent to 
a ruled surface. \qed 
\end{thm}
\noindent
Theorem \ref{thm:rule} follows easily from the development 
of \cite{beau} and, in particular, the classical Noether-Enriques theorem on algebraic surfaces \cite{beau}. 

Two interesting examples of theorem \ref{thm:rule} are the following:
\begin{ex} {\bf ($\exists v\exists y\!\subset\! \exists 
v\forall x\exists y$)}
For any $f\!\in\!\Z[v,y]\!\setminus\!\{0\}$, reconsider its zero set 
as a polynomial in $\Z[v,x,y]$. Abusing notation 
slightly, let us denote this subvariety of $\C^3$ 
by $S_f$. Then there is a natural projection $\varphi$ from 
$S_f$ onto the $v$-axis, and any fiber $\varphi^{-1}(v_0)$ 
is clearly of the form $\{v_0\}\times C$ where $C$ 
is a finite union of lines. So the hypothesis 
of theorem \ref{thm:rule} is satisfied in this example. Better 
still, the conclusion of theorem \ref{thm:rule} is easily verified: 
$S_f$ is clearly birational to a ruled surface, since $S_f$ is 
clearly a Cartesian product of a curve with a line. 
\end{ex} 
\begin{ex}
{\bf ($\exists u\exists v\exists y\!\subset\! \exists u\exists 
v\forall x\exists y$)}
For any $f\!\in\!\Z[u,v,y]\!\setminus\!\{0\}$, reconsider its zero set 
as a polynomial in $\Z[u,v,x,y]$. Abusing notation 
once more, let us denote this subvariety of $\C^4$ 
by $S_f$. Then there is a natural projection 
$\varphi$ from $S_f$ onto the $(u,v)$-plane, and 
any fiber $\varphi^{-1}(u_0,v_0)$ is clearly of the form 
$\{(u_0,v_0)\}\times C$ where $C$ is a finite union 
of lines. More to the point, consider the inverse image  
of $\varphi$ over a line $L$ in the $(u,v)$-plane with positive 
rational slope. Clearly then, $\varphi^{-1}(L)\!\subset\!S_f$ is a 
surface $S'$ fibered over $L$. In particular, every fiber 
$\varphi^{-1}(c)\cap S'$, for $c\!\in\!L$, is clearly a curve 
with a genus zero component. So $S_f$ contains a surface ($S'$) of the 
type specified in Main Theorems 1 and 3, and by theorem 
\ref{thm:rule}, this 
surface has a ruled component.  
\end{ex}
\noindent 

To prove theorem \ref{thm:locusg} we will make some use of 
elimination theory, but in a geometric form. 
\begin{thm}
\label{thm:all} 
Suppose $V$ is a quasi-affine subvariety of 
$\C^m\times \C^n$ defined over $\bQ$, and we respectively use coordinates 
$a\!:=\!(a_1,\ldots,a_m)$ and $x\!:=\!(x_1,\ldots,x_n)$ for the first 
and second factors. Then the following assertions hold: 
\begin{enumerate}
\item{\cite{hart,mumford} The set of $a\!\in\!\C^m$ for which 
there is an $x\!\in\!\Cn$ with $(a,x)\!\in\!V$ is another 
quasi-affine variety $W\!\subseteq\!\C^m$ defined over 
$\bQ$. Furthermore, 
if we are given the sparse encodings of the polynomials 
defining $V$, then we can algorithmically 
determine the analogous data for $W$.} 
\item{\cite{rio,mygcp} Given the sparse encodings of the polynomials defining a 
{\bf zero}-dimensional variety $U\!\subset\!\bQ^n$, the sets 
$U\cap\Zn$ and $U\cap\N^n$ can be effectively computed.}
\item{\cite{gmt,chistov} The decomposition of $V$ into irreducible 
components, and their dimensions, can be effectively computed. \qed} 
\end{enumerate} 
\end{thm} 
\noindent
A more general and explicit version of part (1) 
was derived by Tarski in his work on quantifier 
elimination (over $\R$) in the 1950's. Considerable 
improvements have since been made by other authors, 
giving singly exponential complexity bounds for the 
problem described in part (1). However, since our main concern is 
decidability, we will not dwell on these important 
extensions. 

Finally, we will need Hurwitz' Theorem \cite{miranda,sil} relating the 
genera of the domain and image of a morphism between curves.
\begin{thm}
\label{thm:hurwitz}
Suppose $\varphi : C \longrightarrow C'$ is a nonconstant 
morphism of nonsingular compact curves over $\C$. Let 
$g$ and $g'$ respectively be the genera of $C$ and 
$C'$. Then the following relation holds: 
\[ 2g-2 = (\deg \varphi)(2g'-2)+\sum (e_\varphi(p)-1) \] 
where the sum is over all points $p\!\in\!C$ 
such that $\varphi$ is ramified at $p$, 
and $e_\varphi(p)$ denotes the ramification index. \qed 
\end{thm}

\section{The Proof of Theorem \ref{thm:locusg}}
We will first need the following 
definitions.
\begin{dfn}
Let $\seq(\Zn)$ denote the set of all 
finite multisets of finite subsets of $\Zn$. 
Let us also endow the following partial 
ordering on $\seq(\Zn)$: Declare 
$\{S_1,\ldots,S_k\}\!\leq\!\{T_1,\ldots,T_l\}$ 
iff there are polynomials 
$f_1,\ldots,f_k,g_1,\ldots,g_l\!\in\!\C[x_1,\ldots,x_n]$ 
such that $k\!\leq\!l$, $\supp(f_i)\!=\!S_i$ for all $i$, 
$\supp(g_j)\!=\!T_j$ for all $j$, and 
$\supp(\prod_i f_i)=\supp(\prod_j g_j)$. 
Concluding this connection to factoring, 
let us also define the {\bf factor type}, $\cT_f\!\in\!\seq(\Zn)$, 
of a polynomial $f\!\in\!\C[x_1,\ldots,x_n]$, to be 
the multiset of supports of its irreducible 
factors over $\C[x_1,\ldots,x_n]$. 
\end{dfn} 

It is not hard to see that for a given 
polynomial with parametric coefficients, possessing a 
particular factor type determines a condition 
defining a quasi-affine variety. 
\begin{lemma}
Suppose $f\!\in\!\bQ[a_1,\ldots,a_m,x_1,\ldots,x_n]$. 
Then, for any $\cT\!\in\!\seq(\Zn)$, the set of all 
$(a_1,\ldots,a_m)\!\in\!\C^m$ for 
which $f_a$ has factor type $\cT$ is a quasi-affine variety 
defined over $\bQ$. Furthermore, the polynomials defining this variety are 
effectively computable. 
\end{lemma}
\noindent
{\bf Proof:} First note that $f_a$ has factor type $\geq\!\{S_1,\ldots,S_k\}$ 
iff a set of equations involving $a_1,\ldots,a_m$ 
has a solution. This assertion is immediate, but for 
clarity we give the following example with $k\!=\!2$: 
$a_1+a_2x^2+a_3y^2=(\alpha_1+\beta_1 x+\gamma_1 y)(\alpha_2+
\beta_2 x+\gamma_2 y) \Longleftrightarrow$ the following 
system of equations has a solution $(\alpha_1,\beta_1,\gamma_1,\alpha_2,\beta_2,
\gamma_2)\!\in\!\bQ^6$: 
\begin{eqnarray*} 
\alpha_1 \alpha_2-a_1 & = & 0 \\
\beta_1\beta_2-a_2 & = & 0 \\
\gamma_1\gamma_2-a_3 & = & 0 \\
\alpha_1\beta_2+\alpha_2\beta_1 & = & 0\\
 & \mathrm{etc...} & 
\end{eqnarray*} 

So by theorem \ref{thm:all}, possessing a factor type above or equal to 
$\{S_1,\ldots,S_k\}$ defines a quasi-affine subvariety 
of values of $a$. Now note that the poset of 
possible factor types for any fixed $f$ is finite, 
and recall that quasi-affine varieties are closed under 
any finite sequence of Boolean operations. So by another application of 
theorem \ref{thm:all},  the set of $a$ for which $f_a$ has factor type {\bf exactly} 
$\{S_1,\ldots,S_k\}$ is also a quasi-affine variety. 
So we are done. \qed 

We can now at last prove theorem \ref{thm:locusg}.\\
{\bf Proof of Theorem \ref{thm:locusg}:} Recall once again that 
quasi-affine varieties (even those defined over $\bQ$) are closed under any 
finite sequence of Boolean operations. So it 
suffices to prove that some collection of {\bf algebraic functions} 
of the $a$ in question form a (Turing computable) quasi-affine variety. 
So by lemma 1, it thus suffices to assume that $g$ is an irreducible 
polynomial in $\C[a_1,\ldots,a_m,x,y]$, and $a\!\in\!\bQ^m$ 
is such that $g_a$ is irreducible. 

Let $C$ be the complex zero set of $g_a$. 
Then $C$ is an irreducible curve, possibly with 
singularities. The singularities of $C$  
are precisely the zero set (in the $(x,y)$-plane $\bQ^2$) of an  
effectively constructible system of  
polynomials in $\bQ[a_1,\ldots,a_m,x,y]$ \cite{hart,mumford}. So the 
coordinates of every singular point are algebraic functions of 
$a_1,\ldots,a_m$ defined over $\bQ$. 

Since the number of such singularities is finite, 
a finite sequence of blow-ups will give us a 
new curve $\tC$ birationally equivalent to 
$C$. Furthermore, by our preceding observations, 
and the structure of the blow-up map \cite{hart,mumford}, 
the coefficients of $\tC$ are algebraic functions of 
$a_1,\ldots,a_m$, defined over $\bQ$, as well.

To conclude, note that $\tC$ is a curve in 
$\Pro^N_\C$, where $N$ is $2$ plus the sum of 
the orders of the singularities of $C$. Furthermore, 
we can consider the ambient projective plane 
in which $C$ lies as a coordinate subspace of $\Pro^N_\C$. 
So let $\fii$ be the natural projection mapping 
$\Pro^N_\C$ onto the $x$-axis of this copy of 
$\Pro^2_\C$.

Let us apply theorem \ref{thm:hurwitz} to the preceding morphism 
$\varphi : \tC \longrightarrow \Pro^1_\C$. We then obtain that 
$g_a$ has genus zero iff 
\[ 2\deg \varphi = 2+\sum (e_{\varphi}(p)-1). \] 
Now by theorem \ref{thm:all} and our 
preceding observations (as well as the classical notion 
of discriminant), we can write the preceding sum of ramification 
indices as the order of vanishing of some 
(effectively constructible) polynomial in 
$\bQ[a_1,\ldots,a_m]$ at a point. Furthermore, 
this order of vanishing is equal to some fixed constant iff 
$a$ lies in a quasi-affine variety (defined over $\bQ$) depending on the 
constant. Similarly, the degree of $\varphi$ is some fixed constant iff 
$a$ lies in a quasi-affine variety depending on the constant. 
We thus at last obtain that $\cG_g$ is a Turing constructible 
quasi-affine variety. \qed  

\section{Deciding Prefixes Ending in $\forall\exists$} 
\label{sec:big}

To prove Main Theorems 2 and 3, we will first describe a construction 
which is common to both proofs. 
So let us temporarily consider  
polynomials in $\Z[u,v,x,y]$. More precisely, 
it will be helpful to consider $f$ as a polynomial in $x$ and $y$ with 
coefficients in $\Z[u,v]$. For emphasis, we will 
now respectively write $f_{(u,v)}$ and $f_{(u,v)}(x,y)$ 
in place of $f$ and $f(u,v,x,y)$.  

\begin{dfn}
For any $f\!\in\!\Z[u,v,x,y]$, let $\Xi_f$ be the set of all pairs 
$(u,v)\!\in\!\N^2$ such that 
\[ \forall x \; \exists y 
\mathrm{ \ with \ } f_{(u,v)}(x,y)\!=\!0. \] 
\end{dfn} 

Our main trick for proving Main Theorems 2 and 3 
is the following: create an explicit quasi-affine 
variety $\Omega_f\!\subset\!\C^2$ defined over $\overline{\Q}$, whose 
positive integral points contain (and very nearly equal) 
$\Xi_f$. The following definition and lemma will clarify our 
complex geometric approximation. 
\begin{dfn}
\label{dfn:moreloci} 
Following the notation of definition \ref{dfn:g0}, for any 
$f\!\in\!\C[u,v,x,y]$ let $\Omega_f\!:=\!\cR_f$, 
where we consider $f$ as a polynomial in $x$ and $y$ 
with coefficients in $\C[u,v]$. So $\Omega_f\subset\C^2$ 
and the polynomials defining $\Omega_f$ lie 
in $\C[u,v]$. 
\end{dfn} 
\begin{lemma}
\label{lemma:diop2comp} 
The set $\Xi_f$ is contained in $\Omega_f\cap\N^2$. In particular, 
$\dim\Omega_f\!\geq\!1$ iff $f$ lies in the exceptional locus defined in Main Theorem 3. 

Furthermore, if we restrict to $f\!\in\!\C[v,x,y]$ (and 
thus consider $\Xi_f\!\subseteq\!\N$ and $\Omega_f\!\subseteq\!\C$), 
we have that $\Xi_f\!\subseteq\!\Omega_f\cap\N$. In particular, 
under this restriction, $\dim \Omega_f\!\geq\!1$ iff $f$ lies in the 
exceptional locus defined in Main Theorem 2.
\end{lemma} 

\noindent
{\bf Proof of Lemma \ref{lemma:diop2comp}:} 
By Siegel's Theorem, $(u,v)\!\in\!\Xi_f 
\Longrightarrow S_{f_{(u,v)}}$ contains an irreducible curve of genus zero.  
So the inclusion $\Xi_f\!\subseteq\!\Omega_f\cap\N^2$ is clear. 
The condition for $\dim \Omega_f\!\geq\!1$ is then just a reformulation 
of the condition that enough specializations of $(u,v)$ 
make $S_{f_{(u,v)}}$ have genus zero. The refined statements 
for when $f\!\in\!\Z[v,x,y]$ follow similarly. \qed 

The description of $\Xi_f$ as a subset of the positive 
integral points on a quasi-affine variety allows 
an intuitive complex geometric approach to 
constructing algorithms for a large family of special cases of Diophantine 
prefixes such as $\exists\forall\exists$ and $\exists\exists\forall\exists$. 
In particular, theorem \ref{thm:locusg} tells us that our quasi-affine variety 
$\Omega_f$ is effectively computable, and a judicious use 
of computational algebra can make this implementable.

\begin{rem}
We now briefly clarify the statement of ``genericity'' in 
Main Theorems 2 and 3: Fix the Newton polytope 
$P\!\subset\!\R^3$ of $f$. Then, by theorems  
\ref{thm:kho} and \ref{thm:rule}, and 
the classical Bertini's theorem \cite{hart,mumford}, $S_f$ 
will be an irreducible {\bf non}-ruled surface, provided the coefficients are 
chosen generically and $P$ has at least one lattice point 
in its interior. Thus (except for a meager family of 
supports) Main Theorem 2 implies that for any fixed support we can decide $\exists\forall\exists$ for a Zariski-dense set of $f$. The 
analogous statement for Main Theorem 3 can be derived in 
exactly the same way. 
\end{rem}  

\noindent
{\bf Proof of Main Theorem 2:} 
To construct our necessary algorithm, note that as observed 
in earlier situations, the variable $u$ no longer occurs in $f$. 
So, since the prefixes $\exists u\exists v\forall x\exists y$
and $\exists v\forall x\exists y$ are identical for such
$f$, we may now consider $\Omega_f$ as a subvariety
of $\C$. (And the polynomials defining $\Omega_f$ lie 
in $\bQ[v]$.) Let us also assume $f$ is not identically 
zero. (For when $f$ is identically zero, the prefix in question is trivially 
true.) 

Clearly then, if $\dim \Omega_f\!\leq\!0$, 
deciding whether $\exists v\; \forall x\; \exists 
y$ such that $f(v,x,y)\!=\!0$ reduces to 
simply checking a finite number of instances 
of the prefix $\forall \exists$. (Recall also that 
we can effectively detect $\dim \Omega_f\!\leq\!0$ by theorems \ref{thm:locusg} and \ref{thm:all}.) By the JST 
Theorem and theorem \ref{thm:all}, we thus need only show that the  
hypothesis of Main Theorem 2 implies that 
$\dim \Omega_f\!\leq\!0$. But this follows immediately 
from lemma 2. \qed 

The proof of Main Theorem 3 is almost exactly the 
same, save for the fact that the polynomials 
defining $\Omega_f$ lie in $\bQ[u,v]$. So we will  
omit the proof of Main Theorem 3 and go directly to 
the proof of our first main theorem. 

\section{The Proof of Main Theorem 1} 
Let us temporarily assume that $\biggy_{\N,2}$ is computable. Let 
us also temporarily assume that the stated special 
case of $\exists\exists\forall\exists$ is decidable. 
To derive a contradiction, we will construct an explicit algorithm to 
decide the prefix $\exists\exists\forall\exists$ in 
general. To do this, 
we will again use our algebraic geometric 
trick from the last section. 

Accordingly, our algorithm will have three cases, 
dictated by the topology of $\Omega_f$ in $\C^2$.  
Also note that $\exists u\exists v \forall x\exists y$ 
is trivially true when $f$ is identically zero, so we may assume that $f$ 
is not indentically zero. 

By theorem \ref{thm:locusg}, the JST Theorem, and 
theorem \ref{thm:all} once again, we know that the prefix 
$\forall\exists$ is sufficiently well-behaved so that we can 
make a simplification: we may assume that $\Omega_f$ 
is irreducible. Also, using part (3) of theorem \ref{thm:all}, it is 
clear that Cases I, II, and III below can be distinguished effectively. 
So let us now solve these cases individually. 

\noindent 
{\bf Case I: $\boldsymbol{\Omega_f=\emptyset}$}  

\noindent 
By our preceding observations, we immediately obtain 
that 
\[ \exists u\; \exists v\; \forall x\; \exists y\; 
f(u,v,x,y)\!=\!0\] 
is false. \qed 

\noindent 
{\bf Case II: $\boldsymbol{\dim\Omega_f=0}$}  

\noindent 
Here, we need only check one instance of 
$\forall\exists$. By the JST Theorem, we 
can do this in polynomial time, so 
we are done. \qed 

\noindent
{\bf Case III: $\boldsymbol{\dim\Omega_f\!\geq\!1}$}

\noindent 
By assumption we can solve the case where $\dim\Omega_f\!=\!2$. 
(Indeed, when $\dim\Omega_f\!=\!2$, we can certainly 
find a curve in $\Omega_f$ satisfying the properties 
required in statement (2) of Main Theorem 1.) 
So let us assume $\dim\Omega_f\!=\!1$ and, to simplify notation slightly,   
let $C\!:=\!\Omega_f$. We are left with just two subcases to 
consider and, by assumption, we can compute $\biggy_{\N,2}$ to 
effectively distinguish them. 

{\bf Case III(a): $\boldsymbol{C}$ has finitely many positive integral points}

Since we can compute $\biggy_{\N,2}$, we can simply enumerate 
all possible positive integral points and use the JST Theorem 
a finite (but most likely huge) number of times to 
decide $\exists\exists\forall\exists$. \qed 

{\bf Case III(b): $\boldsymbol{C}$ has infinitely many positive integral 
points} 

By our initial assumption, this case of 
$\exists\exists\forall\exists$ is tractable as well. \qed 

\noindent
Having thus obtained an algorithm contradicting the MR Theorem, 
we are done. \qed 

\section{Conclusion} 
\label{conc}
We have seen a geometric construction which implies 
a weak version of the statement ``$\htp_\N(3)$ 
is decidable $\Longrightarrow \biggy_{\N,2}$ is uncomputable.'' 
The decidability of Hilbert's Tenth Problem 
in three variables is still open, as is 
the existence of computable general upper bounds on 
the size of integral points on algebraic curves. 
So knowing the decidability of $\htp_\N(3)$ or 
the computability of $\biggy_{\N,2}$ would have 
profound implications in algorithmic number theory,  
not to mention arithmetic geometry. 

We emphasize, however, that the uncomputability of $\biggy_{\Z,2}$ 
would by no means contradict the effective upper bounds 
(for heights of integral points) which have already been found 
\cite{bakert,bakerh,bakercoates,brindza,poulaki} 
for certain special classes of curves.  
More precisely, should $\biggy_{\Z,2}$ eventually prove uncomputable, we 
obtain from our development that at least one of the following statements must 
be {\bf true}: 
\begin{itemize}
\item[{\bf (A)}]{ Effective upper bounds on integral points must cease to 
exist for some infinite class of non-superelliptic curves of genus 
at least two.\footnote{Recalling that every curve over $\C$ of genus two 
is hyperelliptic \cite{miranda}, the results of \cite{grant,poon} give 
evidence that this lower bound might need to be increased to three. }}
\item[{\bf (B)}]{ Detecting infinitudes of integral points on a curve 
of genus zero is undecidable.} 
\end{itemize} 
Furthermore, assuming the decidability of detecting infinitudes of 
rational points on curves of genus $\leq\!1$ (and the falsity of 
statement (B)), the uncomputability of $\biggy_{\Z,2}$ would also immediately 
imply the uncomputability of $\biggy_{\Q,2}$. 

It is also clearly the case that the uncomputability of 
$\biggy_{\Z,2}$ would not contradict the decidability 
of $\htp_\Z(2)$, should the latter statement prove true. 
Indeed, the uncomputability of $\biggy_{\Z,2}$ would 
only rule out a stronger version of the decidability of $\htp_\Z(2)$ 
--- the determination of {\bf all} integral points when 
there are only finitely many. More to the point, the existence of effective 
general upper bounds on the height of the {\bf smallest} 
integral point is still an open question. 
For example, Steve Smale has conjectured that 
such upper bounds, for curves of positive genus, 
exist and will be singly exponential in the size of the dense 
encoding \cite{steve}. So the truth of Smale's conjecture would immediately 
imply a brute force algorithm for the positive genus case of $\htp_\Z(2)$.  

We also point out that the exceptional locus in Main Theorem 3  
can be pared down somewhat: Via a suggestion of 
Smale, one can sometimes assume additionally that the forbidden $f$ 
have a zero set which is either (a) reducible or (b) irreducible and singular. 
This refinement is based on examining the critical values 
of (the restriction to $S_f$ of) the natural projection mapping $\C^4$ to 
the $(u,v)$-plane. Other refinements based on a closer 
examination of the real part of $S_f$ are possible and will 
be mentioned in future work. 

We will close by stating a few conjectures and 
open problems related to our development. 
First note that if the statement ``$\biggy_{\N,2}$ 
is computable $\Longleftrightarrow \biggy_{\Z,2}$ is computable'' 
were true, then we could strengthen Main Theorem 1. 
(In particular, we could replace $\biggy_{\N,2}$ by $\biggy_{\Z,2}$.) 
Toward this end, we make the following conjecture on 
a type of equidistribution for integral points on the 
real part of a genus zero curve.
\begin{conj}
Suppose $C\!\subset\!\C^2$ is a curve defined over 
$\Z$ and irreducible over $\C$. Suppose further 
that some irreducible component $C_\R$ of $C\cap\R^2$ has 
noncompact intersection with the 
first quadrant. Then $C_\R$ has infinitely many 
integral points $\Longrightarrow C_\R$ has infinitely 
many positive integral points.
\end{conj} 
\noindent
The truth of this conjecture, combined with a little 
quantifier elimination over $\R$, would immediately 
imply the aforementioned equivalence of $\biggy_{\N,2}$ 
and $\biggy_{\Z,2}$. 

However, a potentially harder problem is to 
refine our proof of Main Theorem 1 to yield 
the truth of the following conjecture we have been 
alluding to.
\begin{conj}
$\htp_\N(3)$ is decidable $\Longrightarrow \biggy_{\N,2}$ 
is uncomputable. 
\end{conj} 

In particular, a refinement of our geometric approach 
seems possible, but quite subtle. For example, 
the proof of a 1970 theorem which essentially computes   
$\biggy_{\Z,2}$ in the special case of genus one curves \cite{bakercoates} 
involves constructing a very special birational map. The map Baker 
and Coates construct takes an arbitrary genus one curve to a 
curve in Weierstrass normal form, preserves rational 
points, and almost preserves integral points. The 
structure of their map was sufficiently good so that 
they could use the previously known bounds for 
curves in Weierstrass normal form to derive height 
bounds for the original (possibly more general) genus one 
curve. 

An analogous construction could be attempted for reducing 
the exceptional locus of Main Theorems 1 and 3 to 
the prefix $\exists\exists\exists$. For instance, 
one could try to use a birational map sending 
$S_f$ to a ruled surface with certain special properties. 
Such a construction, if done properly, could be used to prove 
the undecidability of $\exists\exists\exists$ or the 
equivalence of the decidabilities of $\exists\forall\exists$ 
and $\exists\exists$. Unfortunately, as of 1998, not enough is known 
about integral points on ruled surfaces, or even rational 
surfaces, to make this approach easy. Nevertheless, 
we hope to address this point in the future.  

We also propose the following conjecture motivated 
by our results. 
\begin{conj}
The prefixes $\exists\forall\exists$ and $\exists\exists$ 
are decidable. However, $\biggy_{\N,2}$ is uncomputable 
and $\htp_\N(3)$ is undecidable. 
\end{conj} 
\noindent 
The author is also willing to offer \$1000 (US) for the 
first correct published proof of the decidability of 
$\htp_\N(3)$. This will hopefully prove a safe wager.  

Finally, we remark that for the sake of simplicity, we 
have not given the best possible complexity 
bounds. It is therefore quite likely that Main Theorems 
2 and 3 can be improved to give algorithms which run 
in doubly exponential time. In fact, we are willing 
to conjecture more.  
\begin{conj} 
The Diophantine prefixes $\exists\forall\exists$ and 
$\exists\exists\forall\exists$ are both generically 
decidable within singly exponential time. 
\end{conj} 
\noindent 

\section{Acknowledgements} 
The author would like to express his deep gratitude to 
Professors Lenore Blum, Manuel Blum, and Steve Smale for patiently 
listening to earlier versions of this work. The author 
also thanks James P.\ Jones, Barry Mazur, Zhi-Wei Sun, Shih-Ping Tung, 
and Paul Vojta for useful e-mail discussions. Thanks 
also go to an anonymous referee who made many useful 
suggestions. Special thanks go to Joseph H. Silverman for 
a very encouraging e-mail (see below).  

Finally, the author would like to thank Professor Manuel 
Blum for his extraordinary generosity, and dedicate 
this paper to him. Happy $60$ Manuel! 

\section*{NOTE ADDED IN PROOF}
Joseph H. Silverman has just proved \cite{wow} my Conjecture 1 above, 
so we now have the equivalence of the computabilities of 
$\biggy_{\N,2}$ and $\biggy_{\Z,2}$! We can thus now strengthen Main Theorem 1 
(and sharpen Conjectures 2 and 3) by replacing $\biggy_{\N,2}$ with 
$\biggy_{\Z,2}$ throughout. 

\bibliographystyle{amsalpha}

\begin{thebibliography}{A}

\bibitem[Aya91]{ayad} Ayad, M., {\it 
``On Runge's Theorem,''} Acta.\ Arith.\ {\bf 58} (1991), 
no.\ 2, pp.\ 203--209 (French). 

\bibitem[Bak75]{bakertran} Baker, Alan,  {\it 
Transcendental Number Theory,} Cambridge University Press,  
1975. 

\bibitem[Bak68]{bakert} \underline{\hspace{\baker}},  {\it 
``Contributions to the Theory of Diophantine Equations I:  
On the Representation of Integers by Binary Forms,''} Philos.\ 
Trans.\ Roy.\ Soc.\ London Ser.\ A, 263 (1968), 173--208. 

\bibitem[Bak69]{bakerh} \underline{\hspace{\baker}},  {\it 
``Bounds for the Solutions of the Hyperelliptic Equation,''} 
Proc.\ Camb.\ Philos.\ Soc.\ 65 (1969), 439--444.  

\bibitem[BC70]{bakercoates} Baker, Alan and Coates, John,  {\it 
``Integer Points on Curves of Genus 1,''} Proc.\ Camb. Philos.\ Soc.\ 67 
(1970), 595--602.  

\bibitem[Bea96]{beau} Beauville, Arnaud, {\it 
Complex Algebraic Surfaces,} second edition, London 
Mathematical Society Student Texts, 34, Cambridge 
University Press, 1996, x+132 pp.  

\bibitem[BCSS98]{bcss} Blum, L., Cucker, F., Shub, M., and Smale, 
S.,  {\it Complexity and Real Computation,} foreword by 
Richard M. Karp, Springer-Verlag (1998). 

\bibitem[Bom90]{bomb} Bombieri, Enrico, {\it ``The Mordell Conjecture 
Revisited,''} Ann.\ Sculoa Norm.\ Sup.\ Pisa Cl.\ Sci.\ (4) {\bf 17} 
(1990), no.\ 4, pp.\ 615--640. 

\bibitem[Bri84]{brindza} Brindza, B., {\it ``On $S$-Integral 
Solutions of the Equation $y^m\!=\!f(x)$,''}  
Acta.\ Math.\ Hungar.\ {\bf 44} (1984), no.\ 1--2, pp.\ 
133--139. 

\bibitem[Chi96]{chistov} Chistov, Alexander L., {\it ``Polynomial-Time 
Computation of the Dimension of Algebraic Varieties in 
Zero-Characteristic,''} J.\ Symbolic Comput.\ {\bf 22} 
(1996), no.\ 1, pp.\ 1--25.  

\bibitem[Coh93]{cohen} Cohen, Henri, {\it A Course in Computational  
Number Theory,} Graduate Texts in Mathematics, 138, Springer-Verlag, 
Berlin, 1993, xii+534pp.  

\bibitem[GJ79]{gj} Garey, Michael R.\ and Johnson, David S.\ 
{\it Computers and Intractability: A Guide to the Theory of 
NP-Completeness,} A Series of Books in the Mathematical 
Sciences, W.\ H.\ Freeman and Co., San Francisco, Calif., 
1979, x+338 pp. 

\bibitem[GMT89]{gmt} Gianni, P., Miller, V., and 
Trager, B., {\it ``Decomposition of Algebras,''} 
Symbolic and Algebraic Computation (Rome, 1988), 
pp.\ 300--308, Lecture Notes in Comput.\ Sci., 358, Springer, 
Berlin, 1989. 

\bibitem[Gra94]{grant} Grant, David, {\it ``Integer Points 
on Curves of Genus Two and Their Jacobians,''} Trans.\ 
Amer.\ Math.\ Soc.\ {\bf 344} (1994), no.\ 1, pp.\ 
79--100. 

\bibitem[Har77]{hart} Hartshorne, Robin, {\it Algebraic 
Geometry,''} Graduate Texts in Mathematics, No.\ 52, 
Springer-Verlag.  

\bibitem[Hoe94]{hoeij} van Hoeij, Mark, {\it 
``An Algorithm for Computing an Integral Basis in an 
Algebraic Function Field,''} J.\ Symbolic 
Comput.\ {\bf 18} (1994), no.\ 4, pp.\ 353--363.  

\bibitem[Jon81]{jones81} Jones, James P., {\it ``Classification 
of Quantifier Prefixes Over Diophantine Equations,''}  
Zeitschr.\ f.\ math.\ Logik und Grundlagen d.\ Math., 
Bd.\ 27, pp.\ 403--410 (1981). 

\bibitem[Jon82]{jones9} \underline{\hspace{\jones}}, {\it 
``Universal Diophantine Equation,''} Journal of Symbolic 
Logic, 47 (3), pp.\ 403--410.  

\bibitem[Kho78]{kho78} Khovanskii, A. G.,
{\it ``Newton Polyhedra and the Genus of Complete Intersections,"}
Functional Analysis (translated from Russian), Vol. 12, No. 1,
January--March (1978), pp.\ 51--61.

\bibitem[Lan83]{lang} Lang, Serge, {\it 
Fundamentals of Diophantine Geometry,} Springer-Verlag (1983).  

\bibitem[Mat93]{hilbert10} Matiyasevich, Yuri V., {\it 
Hilbert's Tenth Problem,} foreword by Martin Davis, MIT Press (1993). 

\bibitem[MR74]{matrob} Matiyasevich, Yuri V. and Robinson, Julia {\it 
``Two Universal 3-Quantifier Representations of 
Recursively Enumerable Sets,''} Teoriya Algorifmov i 
Matematicheskaya Logika (Volume dedicated to A. A. Markov), 
pp.\ 112-123, Vychislitel'ny\u{\i} Tsentr, Akademiya Nauk SSSR, Moscow 
(Russian). 

\bibitem[Mir95]{miranda} Miranda, Rick, {\it 
Algebraic Curves and Riemann Surfaces,} Graduate Studies 
in Mathematics, Vol.\ 5, American Mathematical Society.  

\bibitem[Mum95]{mumford} Mumford, David, {\it 
Algebraic Geometry I: Complex Projective Varieties,} 
Reprint of the 1976 edition, Classics in Mathematics, 
Springer-Verlag, Berlin, 1995, x+186 pp.  

\bibitem[Poo96]{poon} Poonen, Bjorn, {\it 
``Computational Aspects of Curves of Genus 
at Least $2$,''} Algorithmic Number Theory (Talence, 1996), 
pp.\ 283--306, Lecture Notes in Comput.\ Sci., 1122, 
Springer, Berlin, 1996. 

\bibitem[Pou93]{poulaki} Poulakis, Dimitrios, {\it 
``Integer Points on Curves of Genus 0,''} Colloq.\ 
Math.\ {\bf 66} (1993), no.\ 1, pp.\ 1--7 (French). 

\bibitem[Roj97]{rio} Rojas, J. Maurice, {\it ``Toric
Laminations, Sparse Generalized Characteristic Polynomials, and a Refinement
of Hilbert's Tenth Problem,''} Foundations of Computational Mathematics,
selected papers of a conference, held at IMPA in Rio de Janeiro, January
1997, Springer-Verlag (1997).

\bibitem[Roj98]{mygcp} \underline{\hspace{\jmr}}, {\it ``Solving Sparse 
Degenerate Polynomial Systems Faster,''} submitted for publication (1998). 
Also available from {\tt http://www.cityu.edu.hk/ma/staff/rojas}. 

\bibitem[Run87]{runge} Runge, C., {\it 
``\"Uber ganzzahlige L\"osungen von Gleichungen zwischen zwei 
Ver\"anderlichen,''} J.\ Reine Angew.\ Math.\ 100 (1887), 
pp.\ 425--435. 

\bibitem[Sch82]{schinzel} Schinzel, Andrzej, {\it 
Selected Topics on Polynomials,} Univ.\ of Michigan 
Press, Ann Arbor, 1982. 

\bibitem[Sch92]{schmidt} Schmidt, Wolfgang M., {\it 
``Integer Points on Curves of Genus 1,''} Compositio 
Mathematica {\bf 81}: 33--59, 1992. 

\bibitem[Sie29]{siegel} Siegel, Carl Ludwig, {\it 
``\"Uber einige Anwendungen Diophantischer Approximationen,''} 
Abh.\ Preuss.\ Akad.\ Wiss.\ Phys.\ Math.\ Kl.\ (1929), Nr.\ 1. 

\bibitem[Sil95]{sil} Silverman, Joseph H., {\it 
The Arithmetic of Elliptic Curves,} corrected 
reprint of the 1986 original, Graduate 
Texts in Mathematics 106, Springer-Verlag (1995).

\bibitem[Sil9\underline{\hspace{\foo}}]{wow} \underline{\hspace{\sil}}, {\it 
``On the Distribution of Integer Points on Curves of 
Genus Zero,''} Theoretical Computer Science, this issue.  

\bibitem[Sma98]{steve} Smale, Steve,  {\it ``Mathematical 
Problems for the Next Century,''} Mathematical 
Intelligencer, to appear (1998).  

\bibitem[Sun92a]{sunphd} Sun, Zhi Wei,  {\it ``Further 
Results on Hilbert's Tenth Problem,''} 
Ph.D. Thesis, Nanjing Univ., Nanjing, 1992.

\bibitem[Tun87]{tungcomplex} Tung, Shih-Ping, {\it 
``Computational Complexities of Diophantine Equations 
with Parameters,''} Journal of Algorithms {\bf 8}, pp.\ 324--336 
(1987).  

\bibitem[Voj87]{vojta} Vojta, Paul, {\it 
Diophantine Approximations and Value Distribution Theory,} 
Lecture Notes in Mathematics, 1239, Springer-Verlag (1987). 

\end{thebibliography}

\end{document}